\documentclass[12pt]{article}

\usepackage{amssymb}
\usepackage{amstex}
\usepackage{amstext}
\usepackage{a4wide}
\usepackage{theorem}
\usepackage[v2,cmtip]{xypic}

\UseComputerModernTips

\numberwithin{equation}{section}

\newtheorem{thm}{Theorem}[section]
\newtheorem{prop}[thm]{Proposition}
\newtheorem{lemma}[thm]{Lemma}
\newtheorem{cor}[thm]{Corollary}

{\theorembodyfont{\rmfamily} \newtheorem{ex}[thm]{Example}}
{\theorembodyfont{\rmfamily} \newtheorem{remark}[thm]{Remark}}

\newcommand{\beginpf}{\smallskip\textbf{Proof. }}
\newcommand{\eopf}{\hfill $\Box$}

\newcommand{\spx}{\ensuremath{\frak{X}}}     
\newcommand{\spy}{\ensuremath{\frak{Y}}}     
\newcommand{\spz}{\ensuremath{\frak{Z}}}     
\newcommand{\spw}{\ensuremath{\frak{W}}}     
\newcommand{\spf}{\ensuremath{\frak{F}}}     

\newcommand{\finite}{\mathcal{F}}      
\newcommand{\sing}{\mathcal{S}}        
\newcommand{\opideal}{\mathcal{I}}     
\newcommand{\allop}{\mathcal{B}}       
\newcommand{\alga}{\mathcal{A}}        
\newcommand{\algb}{\mathcal{B}}        
\newcommand{\algc}{\mathcal{C}}        
\newcommand{\uptri}{\mathcal{T}}       

\newcommand{\nat}{\ensuremath{\Bbb{N}}}        
\newcommand{\inte}{\ensuremath{\Bbb{Z}}}       
\newcommand{\real}{\ensuremath{\Bbb{R}}}       
\newcommand{\complex}{\ensuremath{\Bbb{C}}}    

\newcommand{\sigmaess}{\ensuremath{\operatorname{\sigma_{ess}}}}
\newcommand{\sigmasplit}{\ensuremath{\operatorname{\Sigma_{split}}}}

\newcommand{\gl}{\operatorname{Inv}}
\newcommand{\id}{\operatorname{id}}
\newcommand{\ip}{\operatorname{IP}}
\newcommand{\im}{\operatorname{im}}
\newcommand{\rk}{\operatorname{rk}} 

\newcommand{\ldotscor}{\ensuremath{\ldots@!@!}}

\renewcommand{\phi}{\ensuremath{\varphi}}
\renewcommand{\epsilon}{\ensuremath{\varepsilon}}

\reversemarginpar

\begin{document}
\title{$K$-Theory for Algebras of Operators on Banach Spaces} 
\author{Niels Jakob Laustsen\footnote{Supported by The Danish Natural
Science Research Council.}\\
Department of Mathematics and Computer Science \\
Odense University,
Campusvej 55, DK--5230 Odense M, Denmark\\
 e-mail: \texttt{laustsen@@imada.ou.dk}}
\maketitle

\begin{center}
\textbf{Abstract}
\end{center}
\begin{quote}
We prove that, for each pair $(m,n)$ of non-negative integers,
there is a Banach space $\spx$ for which
$K_0(\allop(\spx))\cong\inte^m$ and $K_1(\allop(\spx))\cong\inte^n.$
Along the way we compute the $K$-groups of all closed ideals of operators
contained in the ideal of strictly singular operators, and we
derive some results about the existence of splittings of certain short exact
sequences. 
\end{quote} 

\section{Introduction and Notation}
Throughout, all vector spaces and algebras are tacitly assumed to be
complex. In addition, all Banach spaces are assumed to be
infinite-dimen\-sional. Banach spaces are denoted by $\spx,\spy,\spz,$
and $\spw$. The term `operator' means a continuous linear map. The
Banach space of all operators from $\spx$ to $\spy$ is denoted by
$\allop(\spx,\spy)$ (or $\allop(\spx)$ in the case where $\spx=\spy$).

By an \emph{ideal of operators} we understand an assignment $\opideal$
which associates to each pair $(\spx,\spy)$ of Banach spaces a subspace
$\opideal(\spx,\spy)$ of $\allop(\spx,\spy)$ satisfying
$$TSR\in\opideal(\spw,\spz)\text{ whenever }
R\in\allop(\spw,\spx),\ S\in\opideal(\spx,\spy), \text{ and }
T\in\allop(\spy,\spz).$$ 
In the case where $\spx = \spy$, we write $\opideal(\spx)$ instead of
$\opideal(\spx,\spx)$. 
An ideal of operators $\opideal$ is said to be \emph{closed} provided that
$\opideal(\spx,\spy)$ is a closed subspace of $\allop(\spx,\spy)$ for
all Banach spaces $\spx$ and $\spy$, and $\opideal$ is said to be {\em
non-zero} if $\opideal(\spx,\spy)\neq\{0\}$ for all (non-zero) Banach
spaces $\spx$ and $\spy$. Let $\opideal_1$ and $\opideal_2$ be two
ideals of operators. We say that $\opideal_1$ is \emph{contained} in
$\opideal_2$ provided that $\opideal_1(\spx,\spy)\subseteq\opideal_2(\spx,\spy)$
for all Banach spaces $\spx$ and $\spy$. Note that an ideal of operators is
non-zero if and only if it contains the ideal $\finite$ of
finite-rank operators. The \emph{rank} of a finite-rank operator $A$ is
denoted by $\rk A$, i.e., $\rk A = \dim(\im A)$.

Recall that an operator $S : \spx\rightarrow\spy$ is said to be
\emph{strictly singular} provided that, whenever $\tilde{\spx}$ is a
subspace of $\spx$ on which $S$ is bounded below, then $\tilde{\spx}$
is finite-dimensional (cf.\ \cite{kato1}). The collection $\sing$ of
strictly singular operators forms a closed ideal of operators
containing the closed ideal of compact operators.

Let $n\in\nat$. We denote the direct sum of the Banach
spaces $\spx_1, \ldotscor , \spx_n$ by
$\bigoplus_{j=1}^n \spx_j$ (or $\spx^n$ in the case where $\spx_1 =  \,\cdots\,  =
\spx_n = \spx$). This is a Banach space under
coordinatewise operations and the maximum norm.

For operators 
$T_j\in\allop(\spx_j,\spy_j)\  (j\in\{1, \ldotscor ,n\})$, we define
$$T_1\oplus \,\cdots\, \oplus T_n : (x_1, \ldotscor ,x_n)\longmapsto
(T_1x_1, \ldotscor ,T_nx_n),\quad \bigoplus_{j=1}^n \spx_j \longrightarrow
\bigoplus_{j=1}^n \spy_j.$$ 
This is clearly a linear and continuous map.

Similar notation and conventions apply to direct sums of (Banach)
algebras and groups.

Let $\alga$ be an algebra, and let $m,n\in\nat$. We denote the vector
space of $(m\times n)$-matrices over $\alga$ by $M_{m,n}(\alga)$ (or
$M_n(\alga)$ in the case where $m=n$). An algebra homomorphism $\phi :
\alga\rightarrow\algb$ induces an algebra homomorphism $\phi_n :
M_n(\alga)\rightarrow M_n(\algb)$ by the definition
$$\phi_n \left(\left(A_{kl}\right)_{k,l=1}^n\right) :=
\left(\phi\left(A_{kl}\right)\right)_{k,l=1}^n.$$

Suppose that $\alga$ is unital. We
write $I$ for the identity of $\alga$, and the group of invertible
elements in $\alga$ is denoted by $\gl(\alga)$. The identity of
$M_n(\alga)$ is denoted by $I^{(n)}$, and we let $\gl_n(\alga) :=
\gl(M_n(\alga))$, the invertible group in $M_n(\alga)$. 

For a non-unital algebra $\alga$, we denote $\alga$ with an
identity adjoined by $\alga^{\sharp}$. Suppose that $\algb$ is a unital
algebra containing $\alga$ as a subalgebra. Then we identify 
$\alga^{\sharp}$ with the subalgebra $\alga + \complex I$ of $\algb$.

Now suppose that $\alga$ is a unital Banach algebra. To make
$M_{m,n}(\alga)$ a Banach space (or a Banach algebra in the case where
$m=n$), we consider it as a closed subspace (subalgebra) of
$\allop(\alga^n,\alga^m)$ in the obvious way.
The norm on $M_{m,n}(\alga)$ thus obtained satisfies
\begin{align*}
\max\big\{\left\| A_{kl}\right\|\,\big|\, k\in\{1,\ldotscor ,m\},\
l\in\{1,\ldotscor ,n\}\big\} &\leq
\left\|\left(A_{kl}\right)_{k,l=1}^{m,n}\right\|\\ 
&\leq \max\bigg\{ \sum_{l=1}^n\left\| A_{kl}\right\|\,\bigg|\,
k\in\{1,\ldotscor ,m\}\bigg\}.
\end{align*}

Let $\spx_1, \ldotscor ,\spx_n$ and $\spy_1, \ldotscor ,\spy_m$ be
Banach spaces. There is a standard bijective correspondance between 
operators
$T\in\allop\left(\bigoplus_{j=1}^n\spx_j,\bigoplus_{j=1}^m\spy_j\right)$
and $(m\times n)$-matrices 
$\left(T_{kl}\right)_{k,l=1}^{m,n}$
with $T_{kl}\in\allop(\spx_l,\spy_k)$, and for every ideal of operators
$\opideal$ we have:
\begin{equation}\label{secondeq}
T\in\opideal\bigg(\bigoplus_{j=1}^n\spx_j,\bigoplus_{j=1}^m\spy_j\bigg)
\text{ if and only if } T_{kl}\in\opideal(\spx_l,\spy_k) \text{ for
all } k,l. 
\end{equation}
In the sequel we \emph{identify} the operator $T$ with its matrix
representation $\left(T_{kl}\right)_{k,l=1}^{m,n}$. In particular,
$\allop(\spx^n,\spx^m)$ is identified with $M_{m,n}(\allop(\spx))$ for
each Banach space $\spx$; under this identification,
$\opideal(\spx^n,\spx^m)$ is identified with
$M_{m,n}(\opideal(\spx))$. 

Note that we have equipped
$M_{m,n}(\allop(\spx))$ with two apparently different norms, one
coming from the identification with $\allop(\spx^n,\spx^m)$ and
another arising from the embedding in
$\allop\left(\allop(\spx)^n,\allop(\spx)^m\right)$. Happily, these
norms coincide, as is easily checked.

\section{Elementary Results}

In this section we prove that both $K$-groups are zero for the algebra
of operators on many of the classical Banach spaces, including $c_0,\,
C(M)\ $($M$ an infinite, compact metric space), $\ell_p$ and
$L_p(\left[0,1\right])\  (p \in \left[1,\infty\right])$. 

First we consider the group $K_0$. Let $\alga$ be an algebra and define
\begin{align*}
\ip(\alga) &:= \{ P \in \alga\,|\, P^2 = P\}, \text{ the set of
idempotents in }\alga, \\ \ip_n(\alga) &:= \ip(M_n(\alga)) \quad (n \in
\nat), \qquad \ip_{\infty}(\alga) := \bigcup_{n \in \nat}
\ip_n(\alga).
\end{align*}
Let $m,n \in \nat$. For $P \in \ip_m(\alga)$ and $Q \in \ip_n(\alga)$
we say that $P \sim_0 Q$ (in $\ip_{\infty}(\alga)$) provided that there 
are matrices $R \in M_{m,n}(\alga)$ and $T \in M_{n,m}(\alga)$
satisfying: $P = RT$ and $Q = TR$. 
Clearly $\sim_0$ is an equivalence relation on
$\ip_{\infty}(\alga)$, and so we may form the quotient $V(\alga) :=
\ip_{\infty}(\alga)/\mathord{\sim_0}$. Let $[P]_V$ denote the 
equivalence class
of $P\in\ip_{\infty}(\alga)$. One easily checks that the operation
$$\left([P]_V,[Q]_V\right)\longmapsto \left[\begin{pmatrix}P & 0 \\ 0
& Q \end{pmatrix}\right]_V,\quad V(\alga)\times
V(\alga)\longrightarrow V(\alga),$$ 
is well-defined and turns $V(\alga)$ into a commutative semigroup. 

Now suppose that the algebra $\alga$ is unital. Then we define
$K_0(\alga)$ to be the Grothen\-dieck group of $V(\alga)$ and, for
$P\in\ip_{\infty}(\alga)$, we denote the canonical image of $[P]_V$ in
$K_0(\alga)$ by $[P]_0$. For $P,Q\in\ip_{\infty}(\alga)$ we observe that
\begin{equation}
[P]_0 = [Q]_0 \text{ if and only if } \begin{pmatrix} P & 0 \\ 0 & I^{(k)}
\end{pmatrix} \sim_0 \begin{pmatrix} Q & 0 \\ 0 & I^{(k)} \end{pmatrix}
\text{ for some } k\in\nat,
\label{thirdeq}
\end{equation}
and we have the following \emph{standard picture} of $K_0(\alga)$:
\begin{equation}
K_0(\alga) = \big\{[P]_0 - [Q]_0\, \big|\, P,Q\in\ip_{\infty}(\alga)\big\}.
\label{sixtheq}
\end{equation}
Moreover, we note that, for each $n\in\nat$ and $P\in\ip_n(\alga)$, the
identity 
\begin{equation}
[P]_0 + [I^{(n)} - P]_0 = [I^{(n)}]_0
\label{fourtheq}
\end{equation}
holds.

In the non-unital case, $K_0(\alga)$ is defined as a subgroup of
$K_0(\alga^{\sharp})$ in the following way. Let $s : \alga^{\sharp}
\rightarrow \alga^{\sharp}$ denote the \emph{scalar map} given by $s(A
+ \zeta I) := \zeta I \  (A \in \alga,\ \zeta \in \complex)$. 
This is clearly an algebra homomorphism, and we define
\begin{equation}
K_0(\alga) := \big\{ [P]_0 - [s_n(P)]_0 \ \big| \ n\in\nat,\, P \in
\ip_n(\alga^{\sharp})\big\},
\label{seventheq}
\end{equation}
where we recall that $s_n$ denotes the map on
$M_n(\alga^{\sharp})$ induced by $s$. 

In the case where $\alga = \allop(\spx)$, the equivalence relation
$\sim_0$ has a nice standard characterization. 

\begin{prop}\label{equi1}
For $P,Q \in \ip_{\infty}(\allop(\spx))$, $P \sim_0 Q$ if and
only if $\im P$ is linearly homeomorphic to $\im Q$.
\eopf
\end{prop}

\begin{lemma}\label{b1}
Suppose that $\spx^n$ is linearly homeomorphic to $\spx^{n + k}$ for some $n,k
\in \nat$. Then
$$K_0(\allop(\spx)) = \big\{ [P]_0 - [Q]_0 \ \big| \ P,Q \in
\ip_n(\allop(\spx))\big\}.$$ 
\end{lemma}
\beginpf 
By \eqref{sixtheq}, it suffices to prove that, for each $m \in \nat$
and $P \in \ip_m(\allop(\spx))$, there is an
operator $\tilde{P} \in \ip_n(\allop(\spx))$ with $P \sim_0
\tilde{P}$. To this end, take $j\in\nat$ for which $n + jk \geq m$. By
assump\-tion, there are operators $R \in M_{n+jk,n}(\allop(\spx))$ and
$T \in M_{n,n+jk}(\allop(\spx))$ satisfying $RT = I^{(n+jk)}$ and $TR =
I^{(n)}$. Let $\tilde{P} := T\,(P\oplus 0)\,R$. Then clearly $\tilde{P} \in
\ip_n(\allop(\spx))$ and $\tilde{P} \sim_0 P\oplus 0\sim_0 P$, as desired.  
\eopf \bigskip 

A Banach space $\spx$ is said to be \emph{primary}
provided that,
whenever $P \in \allop(\spx)$ is idempotent, then $\im P$ or $\im (I - P)$
is linearly homeomorphic to $\spx$.

\begin{prop}\label{b2}
Suppose that $\spx$ is primary and linearly homeomorphic to its square
$\spx^2$. Then $K_0(\allop(\spx)) = \{0\}$. 
\end{prop}
\beginpf
First we note that $[I]_0 = 0$ because, by assumption, $\im I$
is linearly homeomorphic to $\im I^{(2)}$, and so $[I]_0 = 
[I^{(2)}]_0 = [I]_0 + [I]_0$ (cf.\ Proposition~\ref{equi1}). 

By Lemma~\ref{b1}, it suffices
to prove that $[P]_0 = 0$ for each $P \in \ip(\allop(\spx))$. If $\im
P$ is linearly homeomorphic to $\spx$, then $[P]_0 = [I]_0 =
0$. Otherwise, $\im (I - P)$ is linearly homeomorphic to $\spx$, and
the result follows from the identity
\eqref{fourtheq}. 
\eopf 
\begin{remark}\label{gm42a}
As the following example shows, the condition that $\spx$ is linearly
homeomorphic to its square cannot be removed from
Proposition~\ref{b2}, not even if we require $\spx$ to be \emph{prime}
(i.e., the image of \emph{every} idempotent of infinite rank is
linearly homeomorphic to $\spx$), instead of primary.

Let $\spx$ denote the Banach space constructed by Gowers and Maurey in
\cite[\S 4.2]{gm2}. This space is prime and satisfies: for $m,n \in \nat$,
$\spx^m$ is linearly homeomorphic to $\spx^n$ if and only if $m=n$.
It follows immediately from this, Proposition~\ref{equi1}, and
\eqref{thirdeq} that $[I]_0$ is of infinite order in
$K_0(\allop(\spx))$; in particular, $K_0(\allop(\spx)) \neq \{0\}$.

We shall return to this example in Remark~\ref{gm42b}.
\end{remark}

\begin{remark}\label{b2'}
A natural question in relation to Proposition~\ref{b2} is the
following. Suppose that
$\spx$ is linearly homeomorphic to its square and $K_0(\allop(\spx)) =
\{0\}$. Is $\spx$ primary?

The answer to this question is `no'. Take $p,q\in\left[1,\infty\right[$
with $p\neq q$. Then clearly $\ell_p\oplus\ell_q$ is linearly
homeomorphic to its square and $K_0(\allop(\ell_p\oplus\ell_q)) =
\{0\}$ (cf.\ Corollary~\ref{e6}, below), but
$\ell_p\oplus\ell_q$ is not primary, because neither $\ell_p$ nor
$\ell_q$ is linearly homeomorphic to $\ell_p\oplus\ell_q$.
\end{remark} 

We now proceed to consider the group $K_1$. Suppose that $\alga$ is a
unital Banach algebra, and define $\gl_{\infty}(\alga) := \bigcup_{n\in\nat}
\gl_n(\alga)$. Let $m,n\in\nat$. For $U\in\gl_m(\alga)$ and
$V\in\gl_n(\alga)$ we say that $U\sim_1 V$ (in $\gl_{\infty}(\alga)$) 
provided that there are a
natural number $k > \max\{m,n\}$ and a continuous path $t\mapsto
W_t,\, [0,1]\rightarrow \gl_k(\alga)$, with 
$$W_0 = \begin{pmatrix} U & 0 \\ 0 & I^{(k-m)} \end{pmatrix} \qquad \&
\qquad W_1 = \begin{pmatrix} V & 0 \\ 0 & I^{(k-n)} \end{pmatrix}.$$
Clearly $\sim_1$ is an equivalence relation on $\gl_{\infty}(\alga)$,
and so we may define 
$$K_1(\alga) := \gl_{\infty}(\alga)/\mathord{\sim_1}.$$ 
Let $[U]_1$ denote the equivalence class of
$U\in\gl_{\infty}(\alga)$. It is easily checked that the operation
$$\left([U]_1,[V]_1\right)\longmapsto \left[\begin{pmatrix} U & 0 \\ 0
& V \end{pmatrix}\right]_1,\quad K_1(\alga)\times
K_1(\alga)\longrightarrow K_1(\alga),$$
is well-defined and turns $K_1(\alga)$ into a commutative group.

For a non-unital Banach algebra $\alga$, we define $K_1(\alga) :=
K_1(\alga^{\sharp})$. 

\begin{prop}\label{b3}
Suppose that $\spx$ is linearly homeomorphic to its square $\spx^2$ 
and that the
invertible group $\gl(\allop(\spx))$ is connected. Then
$K_1(\allop(\spx)) = \{0\}$.
\end{prop}
\beginpf
Let $n \in \nat$. Clearly $\spx$ is linearly homeomorphic to
$\spx^n$, and so the fact that $\gl(\allop(\spx))$ is connected
implies that $\gl_n(\allop(\spx))$ is connected. Now the result
follows from the definition of $K_1$. 
\eopf

\begin{remark}\label{b4}
In fact Proposition~\ref{b3} is just a special case of a much more
general result. A unital algebra $\alga$ is said to be
\emph{properly infinite} provided that there are elements
$A_1,A_2,B_1,B_2 \in \alga$ satisfying $B_l\,A_k = \delta_{kl}I \ 
(k,l \in \{1,2\})$. One easily checks that the algebra $\allop(\spx)$
is properly infinite in the case where $\spx$ is linearly homeomorphic
to $\spx^n$
for some $n \in \{2,3, \ldotscor \}$. Generalizing the proof of a
$C^*$-algebra result (cf.\ \cite[Lemma 1.2]{cuntz} or \cite[Proposition
3.4]{ror}) yields: for every unital, properly infinite Banach algebra $\alga$,
the group homomorphism $U \mapsto [U]_1,\ \gl(\alga)\rightarrow
K_1(\alga)$, is surjective. In particular, connectedness of $\gl(\alga)$
implies that $K_1(\alga) = \{0\}$.
\end{remark}

\begin{ex}\label{b5}\label{b6}
\begin{enumerate}
\item[(i)] Let $\spx$ denote one of the spaces $c_0$ 
or $\ell_p\  (p \in
\left[1,\infty\right])$. Then certainly $\spx$ is linearly
homeomorphic to its square and, by \cite[Theorems 2.a.3 and 2.a.7]{lt1},
$\spx$ is prime and thus primary. Moreover, the invertible group
$\gl(\allop(\spx))$ is contractible, cf.\ \cite[the Corollaries of
Proposition 2 (p.\ 72) and Lemma 11b (p.\ 79)]{mit}, and consequently
connected. Hence, we obtain
\begin{xxalignat}{2}
& \qquad  K_{\mu}(\allop(c_0)) = \{0\} \qquad \& \qquad
K_{\mu}(\allop(\ell_p)) = \{0\} & (p \in \left[1,\infty\right],\,
\mu\in\{0,1\}).  
\end{xxalignat}
\item[(ii)] Let $p \in \left[1,\infty\right]$. Clearly
$L_p(\left[0,1\right])$ is linearly homeomorphic to its square; it 
is also primary. (This is proved in \cite{m} for $p = 1$ and in
\cite[Theorem 1.3]{aeo} for $1 < p <
\infty$; for $p = \infty$ this follows from (i), because
$L_{\infty}(\left[0,1\right])$ is linearly homeomorphic to
$\ell_{\infty}$, cf.\ \cite{pel}.)\  Furthermore, by \cite[Theorems 2a (p.\
75) and 4 (p.\ 81)]{mit}, $\gl(\allop(L_p(\left[0,1\right])))$ is
contractible and thus connected, so we conclude that
\begin{xxalignat}{2}
&\qquad\qquad\qquad\qquad K_{\mu}(\allop(L_p(\left[0,1\right]))) = \{0\} &
(p\in\left[1,\infty\right],\, \mu\in\{0,1\}).
\end{xxalignat}
\item[(iii)] Let $M$ be an infinite, compact metric space. It is
well-known that $C(M)$ is \mbox{linearly} homeomorphic
to its square. Moreover, $C(M)$ is primary with $\gl(\allop(C(M)))$
contractible. (For $M$ countable, this is proved in \cite[Theorems 1 and
2]{ab}, and for $M$ uncountable, this follows from \cite[Corollary
1]{lp} and \cite[Theorem 2 (p.\ 73)]{mit}.)\  Hence,
\begin{equation} 
K_{\mu}(\allop(C(M))) = \{0\} \tag*{$(\mu\in\{0,1\}).$}
\end{equation}
\end{enumerate}
%
\end{ex}
\bigskip

All the results in this section (except Remark~\ref{gm42a}) are
negative in the sense that they state that the $K$-groups of the
algebra of operators on certain Banach spaces are trivial. To construct
$K$-theoretically more interesting examples, we have to bring the
six-term exact sequence into play. In the next section we prepare
the way for this.

\section{$K$-Theory for the Closed Ideals of Operators Contained in
$\sing$.}\label{singu}

In this section we compute the $K$-groups of $\opideal(\spx)$, where
$\opideal$ is any non-zero, closed
ideal of operators contained in the ideal $\sing$ of strictly
singular operators. Throughout this section, $\opideal$ denotes
such an ideal of operators. The results we obtain
generalize well-known results for the compact operators on a Hilbert
space.

We denote the semigroup of Fredholm operators on $\spx$ by
$\Phi(\spx)$, and, for an operator $T\in\allop(\spx)$, we define the
\emph{essential spectrum} by
$$\sigmaess(T) := \{\zeta\in\complex\ |\ T-\zeta I
\not\in\Phi(\spx)\}.$$
This is a closed, non-empty subset of the spectrum of $T$.

\begin{prop}\label{a2}
Let $n\in\nat$, and let $\left(T_{kl}\right)_{k,l=1}^n\in
M_n(\opideal(\spx)^{\sharp})$. The spectrum of
$\left(T_{kl}\right)_{k,l=1}^n$ is independent of whether it is calculated in
$M_n(\opideal(\spx)^{\sharp})$ or $M_n(\allop(\spx)):$
$$\sigma_{M_n(\opideal(\spx)^{\sharp})}\left(\left(T_{kl}\right)_{k,l=1}^n\right)
=
\sigma_{M_n(\allop(\spx))}\left(\left(T_{kl}\right)_{k,l=1}^n\right),$$
and it is countable.
\end{prop}
\beginpf Write $T_{kl} = S_{kl} + \zeta_{kl}I$, where
$S_{kl}\in\opideal(\spx)$ and $\zeta_{kl}\in\complex$. Since
$\left(S_{kl}\right)_{k,l=1}^n$ is strictly singular, it follows from
\cite[Proposition 2.c.10]{lt1} that
$$\sigmaess\left(\left(T_{kl}\right)_{k,l=1}^n\right) =
\sigmaess\left(\left(\zeta_{kl}I\right)_{k,l=1}^n\right) \subseteq
\sigma_{M_n(\allop(\spx))}\left(\left(\zeta_{kl}I\right)_{k,l=1}^n\right)
\subseteq 
\sigma_{M_n(\complex)}\left(\left(\zeta_{kl}\right)_{k,l=1}^n\right).$$
In particular we see
that $\sigmaess\left(\left(T_{kl}\right)_{k,l=1}^n\right)$ is finite, and 
so, according to \cite[Theorem IV--(5.33)]{kato2},
$\sigma_{M_n(\allop(\spx))}\left(\left(T_{kl}\right)_{k,l=1}^n\right)$
is countable. Now the result follows from
\cite[Corollary (a) of Theorem 10.18]{ru}.
\eopf 
\bigskip 

For a unital Banach algebra $\alga$, we denote the component of
$\gl(\alga)$ containing the identity by $\gl^0(\alga)$.

\begin{lemma}\label{a3}
Let $\alga$ be a unital Banach algebra, and suppose that $A\in\gl(\alga)$ has
countable spectrum. Then $A\in\gl^0(\alga)$.  
\end{lemma}
\beginpf
Take $v\in\left[0,\pi\right[$ with $\sigma(A)\cap(\real
e^{iv}) = \emptyset$. Then $t e^{iv} I +
(1-t)A\in\gl(\alga)$ for all $t\in[0,1]$, and so there is a continuous
path in $\gl(\alga)$ from $A$ to $e^{iv}I\in\gl^0(\alga)$.
\eopf 
\bigskip 

Combining this with Proposition~\ref{a2} and the definition of the
group $K_1$ yields:

\begin{cor}\label{a4}
$K_1(\opideal(\spx)) = \{0\}$. 
\eopf
\end{cor}

Proposition~\ref{a2} has another useful consequence. 
Recall that, for idempotents $P$ and $Q$ in an algebra $\alga$, we
have $P\leq Q$ if
and only if $PQ = P = QP.$ 
The algebra $\alga$ is said to be \emph{finite} provided that,
whenever $P$ and $Q$ are idempotents in $\alga$ satisfying $P\sim_0 Q$
and $P\leq Q$, then $P = Q$. 
If all the algebras $M_n(\alga)\  (n\in\nat)$ are finite, then we say that
$\alga$ is \emph{stably finite}. Note that in the unital case, $\alga$
is finite if and only if every left-invertible (or, equivalently,
right-invertible) element in $\alga$ is invertible.

Now suppose that $\alga$ is a unital Banach algebra which is not finite, and
take $A,B \in \alga$ with $AB = I$ and $BA \neq I$. In particular, $A \neq
0$, and for $\zeta\in\complex$ we have
$$|\zeta| < \frac1{\|A\|} \Longrightarrow I - \zeta A\in\gl(\alga)
\Longleftrightarrow A(B - \zeta I)\in\gl(\alga) \Longrightarrow
\zeta\in\sigma(B).$$
Consequently, $\sigma(B)$ has non-empty interior. Combining this
reasoning with Proposition~\ref{a2} yields:

\begin{prop}\label{sf}
The Banach algebra $\opideal(X)^{\sharp}$ is stably finite.
\eopf
\end{prop}

To compute $K_0(\opideal(\spx))$, we require the following
diagonalization lemma. It is essentially just a restatement of a 
theorem by Edelstein and Wojtaszczyk (cf.\ \cite[Proposition 3.3 and
Theorem 3.5]{ew} or \cite[Theorem
2.c.13]{lt1}), and indeed their proof can be taken over almost
literally; only minor adjustments are necessary to ensure that all the
operators in question belong to (the matrix algebras over)
$\opideal(\spx)^{\sharp}$.

\begin{lemma}\label{a5}
Let $n\in\nat$. For each operator $P\in\ip_n(\opideal(\spx)^{\sharp})$,
there are idempotents $P_1, \ldotscor ,P_n\in\finite(\spx)^{\sharp}$ and an
operator $\Omega\in\gl_n(\opideal(\spx)^{\sharp})$ satisfying
$$\im(\Omega\,P) = \im (P_1 \oplus \,\cdots\,
\oplus P_n).$$
In particular, the equivalence $P\sim_0 P_1 \oplus \,\cdots\, \oplus
P_n$ holds in $\ip_{\infty}(\opideal(\spx)^{\sharp})$. \eopf
\end{lemma}

The following lemma is elementary.

\begin{lemma}\label{equi2}
Let $Q \in \ip_{\infty}(\finite(\spx))$ be of rank one, and let $n\in\nat$.
For every $P \in \ip_n (\finite(\spx))$, the formulas
$$[P]_0 = (\rk P)\,[Q]_0 \qquad \& \qquad [I^{(n)} - P]_0 - [I^{(n)}]_0 = -(\rk
P)\,[Q]_0$$ 
hold in $K_0(\opideal(\spx)^{\sharp})$. In particular,
$[P]_0 - [s_n(P)]_0 \in \inte\,[Q]_0$ for every
$P\in\ip(M_n(\finite(\spx))^{\sharp})$. 
\mbox{  }\eopf 
\end{lemma}

\begin{thm}\label{k0s}
Suppose that $Q \in \ip_{\infty}(\finite(\spx))$ is of rank one. Then the map 
$$\omega : \nu \longmapsto \nu\,[Q]_0,\quad \inte \longrightarrow
K_0(\opideal(\spx)),$$
is a group isomorphism.
\end{thm}  
\beginpf
Clearly $\omega$ is additive.

Assume towards a contradiction that $\omega$ is not injective. Then
there is a natural number $\nu$ satisfying  
$$0 = \nu\,[Q]_0 = [\underbrace{Q \oplus \,\cdots\, \oplus
Q}_{\nu}]_0$$
in $K_0(\opideal(\spx))$ and thus in $K_0(\opideal(\spx)^{\sharp})$,
and consequently $0 \oplus \,\cdots\, \oplus 0 \oplus I^{(k)} 
\sim_0 Q \oplus \,\cdots\,
\oplus Q \oplus I^{(k)}$ for some $k \in \nat$, cf.\
\eqref{thirdeq}. Clearly $0 \oplus \,\cdots\,
\oplus 0 \oplus I^{(k)} \leq Q \oplus \,\cdots\, \oplus Q \oplus
I^{(k)}$, and so, by
Proposition~\ref{sf}, $0 \oplus \,\cdots\, \oplus 0 \oplus I^{(k)} = Q \oplus
\,\cdots\, \oplus Q \oplus I^{(k)}$, contradicting the fact that $Q \neq 0$. 
Hence, $\omega$ is injective. 

To show that $\omega$ is surjective, let $g \in K_0(\opideal(\spx))$
be given. By \eqref{seventheq}, we can take $n \in \nat$ and $P \in
\ip_n(\opideal(\spx)^{\sharp})$ so that $g = [P]_0 - [s_n(P)]_0$, and,
by Lemma~\ref{a5}, we can take $P_1, \ldotscor ,P_n \in
\ip(\finite(\spx)^{\sharp})$ for which $P \sim_0 P_1 \oplus \,\cdots\,
\oplus P_n$ in $\ip_{\infty}(\opideal(\spx)^{\sharp})$. Consequently, 
$$g = [P_1 \oplus \,\cdots\, \oplus P_n]_0 - [s_n(P_1 \oplus \,\cdots\,
\oplus P_n)]_0  = \sum_{j=1}^n \left([P_j]_0 - [s(P_j)]_0\right)
\in \inte\, [Q]_0,$$ 
cf.\ Lemma~\ref{equi2}.
\eopf 

\section{Splittings and H.I.\ Spaces} 
\label{splitting}
Let $\alga$, $\algb$, and $\algc$ be algebras, and suppose that we have
a short exact sequence
$$\diagram
\Sigma : \qquad \{0\} \rto &\alga \rto^{\displaystyle{\phi}} &\algb
\rto^-{\displaystyle{\psi}} &\algc \rto &\{0\}.
\enddiagram$$
We say that $\Sigma$ \emph{splits algebraically} (or is \emph{split exact})
provided that there is an algebra homomorphism $\theta :
\algc\rightarrow \algb$ for which $\psi\circ\theta = \id_{\algc}$;
such a map $\theta$ is called an \emph{(algebraic) splitting
homomorphism}. The functor $K_0$ does not in general preserve short
exact sequences, but it does preserve split exact sequences.

Now suppose that $\alga$, $\algb$, and $\algc$ are Banach algebras and
that $\phi$ and $\psi$ are continuous. 
We say that $\Sigma$ \emph{splits strongly} (or is
\emph{strongly split exact}) provided that there is a con\-tinuous
splitting homomorphism. The functor $K_1$ does not in general preserve
short exact sequences, but it does preserve strongly split exact sequences.

A fundamental theorem of $K$-theory states that, to every short exact
sequence $\Sigma$, where $\alga$, $\algb$, and $\algc$ are Banach
algebras and $\phi$ and $\psi$ are continuous, we
can associate the cyclic six-term exact sequence
$$\spreaddiagramrows{2pc}\spreaddiagramcolumns{2pc}
\diagram
\displaystyle{K_0(\alga)}\rto^-{\displaystyle{K_0(\phi)}} & 
\displaystyle{K_0(\algb)}\rto^-{\displaystyle{K_0(\psi)}} &
\displaystyle{K_0(\algc)}\dto^-{\displaystyle{\delta_0}} \\
\displaystyle{K_1(\algc)}\uto^-{\displaystyle{\delta_1}} &
\displaystyle{K_1(\algb)}\lto_-{\displaystyle{K_1(\psi)}} &
\displaystyle{K_1(\alga).}\lto_-{\displaystyle{K_1(\phi)}}
\enddiagram$$
The map $\delta_1$ is called the \emph{index map} because of its
relation to the Fredholm index $i$, described in Proposition \ref{d3},
below.

In this section we shall study the conditions under which the short
exact sequence
$$\diagram
\Sigma_{\spx} : \qquad \{0\} \rto &\opideal(\spx) 
\rto^{\displaystyle{\iota}} &\allop(\spx)
\rto^-{\displaystyle{\pi}} &\allop(\spx)/\opideal(\spx) \rto &\{0\}
\enddiagram$$
splits. Here, as in the rest of this section, $\opideal$ denotes a
closed, non-zero ideal of operators contained in $\sing$. Note that in
this case, by the results of \S\ref{singu}, the six-term exact
sequence has the following form
$$\spreaddiagramcolumns{1pc}\spreaddiagramrows{-1pc}
\diagram
\displaystyle{\{0\}}\rto &
\displaystyle{K_1(\allop(\spx))}\rto^-{\displaystyle{K_1(\pi)}} &
\displaystyle{K_1(\allop(\spx)/\opideal(\spx))}
\rto^-{\displaystyle{\delta_1}}&
\displaystyle{K_0(\opideal(\spx)) = \inte\,[Q]_0}
\rto^-{\displaystyle{K_0(\iota)}} & \cdots \\
& \qquad \cdots\rto^-{\displaystyle{K_0(\iota)}}&
\displaystyle{K_0(\allop(\spx))}
\rto^-{\displaystyle{K_0(\pi)}}& 
\displaystyle{K_0(\allop(\spx)/\opideal(\spx))}\rto& 
\displaystyle{\{0\}},
\enddiagram$$
where $Q$ is an idempotent of rank one.

Our first result is well-known for Hilbert spaces, and the proof of
this special case is easily adapted to cover the Banach-space case.

\begin{prop}\label{d3}
For every $m\in\nat$, the diagram
$$\spreaddiagramcolumns{1.5pc}\spreaddiagramrows{2pc}
  \diagram
  \displaystyle{\Phi(\spx^m)} \rto|>>\tip^-{\displaystyle{\pi_m}} 
  \dto_{\displaystyle{i}} &
  \displaystyle{\gl_m(\allop(\spx)/\opideal(\spx))}
  \rto^-{\displaystyle{[\,\cdot\,]_1}} &
  \displaystyle{K_1(\allop(\spx)/\opideal(\spx))}
  \dto^{\displaystyle{\delta_1}} \\
  \displaystyle{\inte}\rrto^{\displaystyle{\omega}}_{\displaystyle{\cong}} &
  &\displaystyle{K_0(\opideal(\spx))}
\enddiagram$$
commutes.
\eopf
\end{prop}

As a consequence, we obtain:
\begin{prop}\label{d4}
The following assertions are equivalent:
\begin{enumerate}
\item[\emph{(a)}] the map $K_0(\iota) : K_0(\opideal(\spx))\longrightarrow
K_0(\allop(\spx))$ is injective;
\item[\emph{(b)}] the map $\delta_1 : K_1(\allop(\spx)/\opideal(\spx))
\longrightarrow K_0(\opideal(\spx))$ is zero;
\item[\emph{(c)}] the map $K_1(\pi) : K_1(\allop(\spx))\longrightarrow
K_1(\allop(\spx)/\opideal(\spx))$ is surjective (and thus an
isomorphism);
\item[\emph{(d)}] for each $m\in\nat$, $\spx^m$ is not linearly homeomorphic
to any of its proper, closed subspaces of finite codimension;
\item[\emph{(e)}] for each $m\in\nat$, every Fredholm operator on $\spx^m$
has Fredholm index zero.
\end{enumerate}
These assertions are fulfilled in the case where the algebra $\allop(\spx)$ is
stably finite.
\end{prop}
\beginpf
The fact that (a), (b), and (c) are equivalent follows directly from the
six-term exact sequence; note that, by Corollary~\ref{a4}, $K_1(\pi)$
is always injective.

The implication `(b) $\Rightarrow$ (e)' is immediate from
Proposition~\ref{d3}.

We prove that (d) implies (a) by contraposition. Suppose that
$K_0(\iota)$ is not injective. Then there exist a number $j\in\nat$
and an idempotent $P\in\finite(\spx^j)\setminus\{0\}$ for which $[P]_0
= 0$ in $K_0(\allop(\spx))$. By \eqref{thirdeq} and \eqref{fourtheq}, this
implies that $(I^{(j)} - P)\oplus I^{(k)} \sim_0 I^{(j)}\oplus
I^{(k)}$ in $\ip_{\infty}(\allop(\spx))$ for some
$k\in\nat$ so, by Proposition~\ref{equi1}, $\im((I^{(j)} - P)\oplus 
I^{(k)})$ is
linearly homeomorphic to $\im(I^{(j)}\oplus I^{(k)}) = \spx^{j+k}$. Clearly
$\im((I^{(j)} - P)\oplus I^{(k)})$ is a proper, closed subspace of codimension
$\rk P$ in $\spx^{j+k}$.

The implication `(e) $\Rightarrow$ (d)' is also proved
contrapositively. Suppose that there exist a number $m\in\nat$ and a
proper, closed subspace $\spw$ of finite codimension in $\spx^m$ which
is linearly homeomorphic to $\spx^m$. Take a linear homeomorphism
$\tilde{T} : \spw\rightarrow \spx^m$, and let $\spf$ be a
finite-dimensional complement of $\spw$ in $\spx^m$. Then $T :=
\tilde{T}\oplus\ 0|_{\spf} \in \allop(\spx^m)$ is a Fredholm operator
of index $i(T) = \dim\spf\geq 1$, and so (e) is not satisfied. Moreover,
the last remark follows from this, for if $R := \tilde{T}^{-1} :
\spx^m\rightarrow \spw$ is considered as an operator on $\spx^m$, then
$TR = I^{(m)}$, but $RT\neq I^{(m)}$, so the algebra $\allop(\spx^m)$ is
not finite.
\eopf

\begin{cor}\label{d5}
\label{split-exact}
\begin{enumerate}
\item[\emph{(i)}] Suppose that the short exact sequence 
$\Sigma_{\spx}$ splits
algebraically. Then $K_0(\allop(\spx))$ contains a subgroup isomorphic
to $\inte$. In particular, $\Sigma_{\spx}$ never splits algebraically in the
case where $K_0(\allop(\spx)) = \{0\}$.
\item[\emph{(ii)}] Suppose that there exists a number $m\in\nat$ for which
$\spx^m$ is linearly homeomorphic to one of its proper, closed
subspaces of finite codimension, or, equivalently, admits a Fredholm
operator of non-zero index. Then $\Sigma_{\spx}$ does not split algebraically. 
\eopf
\end{enumerate}
\end{cor}
 
\begin{remark}\label{gm42b}
It follows from (ii) that the condition in (i) is not
sufficient for an algebraic splitting of the short exact
sequence $\Sigma_{\spx}$ to
exist. Indeed, consider the prime Banach space $\spx$ introduced in
Remark~\ref{gm42a}. We saw there that the map
$$\nu \longmapsto \nu\,[I]_0,\quad \inte\longrightarrow K_0(\allop(\spx)),$$
is a group
monomorphism, but $\Sigma_{\spx}$ does not split algebraically, because the
primeness of $\spx$ implies that $\spx$ is linearly homeomorphic to
each of its subspaces of finite codimension. 

We do not know whether the equivalent conditions in
Proposition~\ref{d4} are sufficient for an algebraic splitting of
$\Sigma_{\spx}$ to exist, but it seems unlikely.
\end{remark}
\medskip

The result of Corollary~\ref{d5} might give the impression that the
short exact sequence $\Sigma_{\spx}$ \emph{never} splits. This is,
however, not true.  Recall that a Banach space $\spx$ is said to be
\emph{hereditarily indecomposable} or an \emph{H.I.\ space} provided
that no closed, infinite-dimensional subspace $\tilde{\spx}$ of $\spx$
admits an idempotent $P\in\allop(\tilde{\spx})$ for which neither $P$
nor $I|\tilde{\spx} - P$ is of finite rank; an equivalent definition
is that, for each $c>0$ and each pair $(\spx_1,\spx_2)$ of
infinite-dimensional subspaces of $\spx$, there are unit vectors
$x_1\in\spx_1$ and $x_2\in\spx_2$ with $\| x_1 - x_2 \| \leq c$.
Gowers and Maurey's fundamental results about H.I.\ spaces are that
they exist and that, for every H.I.\ space $\spx$, $\allop(\spx) =
\sing(\spx) + \complex I$ (cf.\ \cite{gm1}). In particular, this
implies that $\Sigma_{\spx}$ splits strongly for $\opideal =
\sing$. We shall now extend this result and derive its $K$-theoretical
consequences.

For $m\in\nat$, we denote the closed subalgebra of $M_m(\complex)$ of
upper triangular matrices by $\uptri_m(\complex)$.

\begin{prop}\label{s1}
For each $m\in\nat$, there exist a Banach space $\spy$ and continuous algebra
homomorphisms $\psi$ and $\theta$ giving a strongly split exact sequence
$$\diagram
\sigmasplit : \qquad 
\{0\} \rto &\sing(\spy) \rto^{\displaystyle{\iota}} &\allop(\spy)\rto<0.33ex>^{\displaystyle{\psi}}
&\uptri_m(\complex) \lto<0.33ex>^{\displaystyle{\theta}} \rto &\{0\}.
\enddiagram$$
\end{prop}
\beginpf
Let $\spy_1$ be an H.I.\ space, and take a descending chain
$\spy_1\supseteq\spy_2\supseteq \,\cdots\, \supseteq\spy_m$ of closed
subspaces for which $\spy_{j+1}$ is uncomplemented in $\spy_j\ 
(j\in\{1, \ldotscor ,m-1\})$. (The existence of a closed, uncomplemented subspace
$\spy_{j+1}$ in $\spy_j$ follows from \cite[Theorem 1]{lt3}, because
$\spy_j$ inherits the H.I.\ property from $\spy_1$ and is
consequently not linearly homeomorphic to a Hilbert space.)\ 
Let $\spy := \bigoplus_{j=1}^m \spy_j$, and consider an operator $T =
\left(T_{kl}\right)_{k,l=1}^m \in \allop(\spy)$.

For $1\leq k\leq l\leq m$, $\spy_l$ is a closed subspace of the H.I.\ 
space $\spy_k$ so, by \cite{fe}, every operator $T_{kl} \in
\allop(\spy_l,\spy_k)$ is of the form $T_{kl} = S_{kl} + \zeta_{kl}
I_{kl}$, where $S_{kl} \in \sing(\spy_l,\spy_k),\,
\zeta_{kl} \in \complex$, and $I_{kl}$ denotes the inclusion
map of $\spy_l$ into $\spy_k$; in particular, $I_{kk} = I|\spy_k$.

For $1\leq l<k\leq m$, we have $\spy_k \subseteq \spy_l$, and so we 
may consider
$T_{kl}$ as an operator on $\spy_l$ with $\im T_{kl} \subseteq
\spy_k$. Hence, $T_{kl} = S_{kl} + \zeta_{kl} I|\spy_l$, where $S_{kl}$
is strictly singular and $\zeta_{kl} \in \complex$. Since $\spy_k$ is
of infinite codimension in $\spy_l$, $\im T_{kl}$ must be of infinite
codimension in $\spy_l$ as well, so $\zeta_{kl} = 0$, i.e., $T_{kl} =
S_{kl} \in \sing(\spy_l,\spy_k)$.

Hence, we conclude that every operator on $\spy$ is of the form
\begin{equation*}
T = \begin{pmatrix} 
  \zeta_{11} I + S_{11} & \zeta_{12} I_{12} + S_{12}          & 
                 \ldots & \zeta_{1m} I_{1m} + S_{1m}          \\
                 S_{21} & \zeta_{22} I|\spy_2 + S_{22}             &
                 \ldots & \zeta_{2m} I_{2m} + S_{2m} \\
                 \vdots & \vdots    &    \ddots & \vdots           \\
                 S_{m1} &   S_{m2}  & \ldots                       & 
                 \zeta_{mm} I|\spy_m + S_{mm} \end{pmatrix},
\end{equation*}
where $\zeta_{kl} \in \complex$ and $S_{kl} \in
\sing(\spy_l,\spy_k)$. The complex numbers 
$\zeta_{kl}$ are uniquely determined by $T$. Consequently, it makes
sense to define a map 
$$\psi : T \longmapsto
  \begin{pmatrix} 
  \zeta_{11} & \zeta_{12} & \ldots & \zeta_{1m} \\
       0     & \zeta_{22} & \ldots & \zeta_{2m} \\
    \vdots   &    \vdots  & \ddots & \vdots     \\
       0     &     0      & \ldots & \zeta_{mm} \end{pmatrix},\quad
  \allop(\spy)\longrightarrow\uptri_m(\complex).$$
Direct calculations show that $\psi$ is an
algebra homomorphism. By \eqref{secondeq}, $\ker\psi = \sing(\spy)$.

Since the map
$$\theta : \begin{pmatrix} 
  \zeta_{11} & \zeta_{12} & \ldots & \zeta_{1m} \\
       0     & \zeta_{22} & \ldots & \zeta_{2m} \\
    \vdots   &   \vdots   & \ddots & \vdots     \\
       0     &     0      & \ldots & \zeta_{mm} \end{pmatrix} \longmapsto
  \begin{pmatrix} 
  \zeta_{11}I & \zeta_{12}I_{12} & \ldots & \zeta_{1m}I_{1m} \\
        0     & \zeta_{22}I|\spy_2    & \ldots & \zeta_{2m}I_{2m} \\
    \vdots   &   \vdots   & \ddots & \vdots     \\
       0     &     0      & \ldots & \zeta_{mm}I|\spy_m
  \end{pmatrix},\quad
  \uptri_m(\complex) \longrightarrow \allop(\spy),$$
is an algebra homomorphism for which $\psi\circ\theta =
\id_{\uptri_m(\complex)}$, we have completed the construction of the
split exact sequence $\sigmasplit$; $\theta$ is continuous because
$\uptri_m (\complex)$ is finite-dimensional, and $\psi$ is continuous
because of the following easy result.
\eopf

\begin{lemma}\label{s2}
Let $\algb$ and $\algc$ be Banach algebras, let $\opideal$ be a closed
ideal in $\algb$, and suppose that we have a short exact sequence
$$\diagram
\{0\} \rto &\opideal \rto^{\displaystyle{\iota}} &\algb \rto^{\displaystyle{\psi}} &\algc \rto &\{0\}.
\enddiagram$$
Then $\algc$ is isomorphic to the quotient $\algb/\opideal$ via the
algebra isomorphism $\chi : B + \opideal \mapsto \psi(B),\,
\algb/\opideal\rightarrow \algc$. Moreover, $\psi$ is continuous if and
only if $\chi$ is continuous.
\eopf
\end{lemma}

Since $\uptri_m(\complex)$ is homotopy equivalent to its diagonal
subalgebra $\complex^m$, these two algebras have isomorphic
$K$-groups, i.e., $K_0(\uptri_m(\complex))\cong\inte^m$ and
$K_1(\uptri_m(\complex)) = \{0\}$.

This, together with the strongly split exact sequence $\sigmasplit$ and
strong split exactness of the $K$-functors, gives:
\begin{cor}\label{s6}
Let $\spy$ be the above Banach space. Then
$$K_0(\allop(\spy))\cong\inte^{m+1} \qquad \& \qquad 
K_1(\allop(\spy)) = \{0\}.$$
The generators of $K_0(\allop(\spy))$ are $[Q]_0, [P_1]_0,
\ldotscor ,[P_m]_0$, where $Q$ is an idempotent of rank one and $P_j$
denotes the idempotent operator on $\spy$ given by the $(m\times
m)$-matrix with $I|\spy_j$ in position $(j,j)$ and zeroes everywhere
else.
\eopf 
\end{cor}

\section{The Main Result}\label{e}

Recall that
Banach spaces $\spx$ and $\spy$ are said to be \emph{totally
incomparable} provided that no infinite-dimensional subspaces $\tilde{\spx}$
and $\tilde{\spy}$ of $\spx$ and $\spy$, respectively, are linearly
homeomorphic (cf.\ \cite{ro}).
\begin{ex}\label{e4}
\begin{enumerate}
\item[(i)] By \cite[p.\ 75]{lt1}, any two distinct spaces in the
family $\{\ell_p\, |\, p\in\left[1,\infty\right[\}\cup \{c_0\}$ are
totally incomparable.
\item[(ii)] Suppose that $\spx$ has an unconditional basis and $\spy$ is an
H.I.\ space. Then $\spx$ and $\spy$ are totally incomparable. This
follows from the fact that every closed, infinite-dimensional subspace
$\tilde{\spx}$ of $\spx$ contains an unconditional basic sequence 
(cf. \cite[the comment just after Problem 1.d.5]{lt1}), whereas $\spy$
contains no unconditional basic sequence (cf.\ \cite{gm1}).
\end{enumerate}
\end{ex}
We are now ready to state and prove our main theorem. For convenience,
we use the convention that $\inte^0 := \{0\}$, the trivial group.
\begin{thm}\label{e5}
For every pair $(m,n)$ of non-negative integers, there is a Banach
space $\spx$ for which 
$$K_0(\allop(\spx))\cong\inte^m \qquad \& \qquad
K_1(\allop(\spx))\cong\inte^n.$$
\end{thm}
\beginpf
Take $\spx_1, \ldotscor ,\spx_{n+1}$ to be distinct spaces from the family
$\{\ell_p\, |\, p\in\left[1,\infty\right[\} \cup \{c_0\}$. It follows
immediately from the six-term exact sequence and Example~\ref{b6} (i)
that 
\begin{xxalignat}{2}
&\quad \  
K_0\left(\allop(\spx_j)/\sing(\spx_j)\right) =
\{0\} \qquad \& \qquad K_1\left(\allop(\spx_j)/\sing(\spx_j)\right)
\cong \inte & \left(j\in\{1, \ldotscor ,n + 1\}\right).
\end{xxalignat}
If $m = 0$, let $\spx_{n+2} := \{0\}$. Otherwise, let
$\spx_{n+2} := \spy$, the space constructed in the proof of 
Proposition~\ref{s1}. Note that in both cases
$$K_0\left(\allop(\spx_{n+2})/\sing(\spx_{n+2})\right) \cong
\inte^m \qquad \& \qquad 
K_1\left(\allop(\spx_{n+2})/\sing(\spx_{n+2})\right) = \{0\}.$$

Let $\spx := \bigoplus_{j=1}^{n+2}\spx_j$. By
Example~\ref{e4} (i)--(ii) and \cite[Corollary 4 (b)]{ro}, $\spx_1, \ldotscor
,\spx_{n+2}$ are mutually totally incomparable, and so
$\allop(\spx_k,\spx_l) = \sing(\spx_k,\spx_l)\  (k\neq l)$. 
Consequently, the obvious diagonal embedding of
$\bigoplus_{j=1}^{n+2}\allop(\spx_j)/\sing(\spx_j)$ into
$\allop(\spx)/\sing(\spx)$ is in fact a homeomorphic algebra
isomorphism; in particular, this implies that
$$K_{\mu}\left(\allop(\spx)/\sing(\spx)\right) \cong
\bigoplus_{j=1}^{n+2}K_{\mu}\left(\allop(\spx_j)/\sing(\spx_j)\right)
\cong \begin{cases} \inte^m& \text{for } \mu = 0,\\
                    \inte^{n+1}& \text{for } \mu = 1. \end{cases}$$
Combining this with the six-term exact sequence gives the following diagram
$$\spreaddiagramrows{2pc}\spreaddiagramcolumns{1.2pc}
\diagram
\displaystyle {\inte} \rto^-{\displaystyle{\omega}}_-{\displaystyle{\cong}} 
&\displaystyle{K_0(\sing(\spx))} \rto^-{\displaystyle{K_0(\iota)}} 
&\displaystyle{K_0(\allop(\spx))} \rto^-{\displaystyle{K_0(\pi)}}
&\displaystyle{K_0(\allop(\spx)/\sing(\spx))}
\dto & \displaystyle{\hspace{-8ex} \cong \ \inte^m} \\
\displaystyle{\inte^{n+1}}\uto^-{\displaystyle{\sigma}}
&\displaystyle{K_1(\allop(\spx)/\sing(\spx))}
\lto^-{\displaystyle{\cong}}\uto_-{\displaystyle{\delta_1}} & 
\displaystyle{K_1(\allop(\spx))}
\lto_-{\displaystyle{K_1(\pi)}} & \{0\}, \lto
\enddiagram$$
where the right-hand-side rectangle is exact, and $\sigma$ denotes the
induced group homomorphism, so that the left-hand-side square
commutes. In fact $\sigma$ is given by addition: $\sigma
(\nu_1,\ldotscor ,\nu_{n+1}) = \sum_{j=1}^{n+1}\nu_j$. This follows
from a straightforward calculation based on the fact that, for
Fredholm operators $T_j\in\allop(\spx_j)\  (j\in\{1,\ldotscor,n +
2\})$, the operator $T_1\oplus\,\cdots\,\oplus T_{n+2} \in
\allop(\spx)$ is a Fredholm operator of index $\sum_{j=1}^{n+2}
i(T_j)$.  In particular, $\sigma$ is surjective, so, by exactness,
$K_0(\pi)$ is injective and thus an isomorphism, i.e.,
$K_0(\allop(\spx)) \cong \inte^m$.

Moreover, we have that $K_1(\allop(\spx)) \cong \im K_1(\pi) =
\ker\delta_1 \cong \ker\sigma \cong \inte^n$.
\eopf 
\bigskip 

Careful tracking of the isomorphisms in the above proof shows that the
generators of $K_0(\allop(\spx))$ and $K_1(\allop(\spx))$ are given as
follows. 

The group $K_0(\allop(\spx))$ is generated by the elements
$[\underbrace{0\oplus\,\cdots\,\oplus 0}_{n+1}\oplus P_j]_0\  (j\in\{1,
\ldotscor ,m\})$,
where the idempotents $P_j$ are defined as in Corollary~\ref{s6}.

The generators of $K_1(\allop(\spx))$ are $[U_j]_1\  (j\in\{1, \ldotscor
,n\})$ with the invertible operator $U_j$ given by the matrix
$$\begin{pmatrix}
I &        &   &     &   &        &   &   &   \\
  & \ddots &   &     &   &        &   &   &   \\
  &        & I &     &   &        &   &   &   \\
  &        &   & L_j &   &        &   &   &   \\
  &        &   &     & I &        &   &   &   \\
  &        &   &     &   & \ddots &   &   &   \\
  &        &   &     &   &        & I &   &   \\
  &        &   & A_j &   &        &   & R &   \\
  &        &   &     &   &        &   &   & I\end{pmatrix}.$$
Here $L_j$ (in position $(j,j)$) denotes the unilateral left-shift on
$\spx_j$, $R$ (in position $(n + 1,n + 1)$) denotes the unilateral
right-shift on $\spx_{n+1}$, and $A_j$ (in position $(n + 1,j)$) is the
rank-one operator defined by $(\zeta_k)_{k\in\nat}\mapsto
(\zeta_1,0,0, \ldotscor ),\, \spx_j\rightarrow\spx_{n+1}$.
\bigskip

As a special case of Theorem~\ref{e5} and its proof we note that:
\begin{cor}\label{e6}
Let $n\in\nat$, and let $\spx_1, \ldotscor ,\spx_{n+1}$ be distinct
spaces belonging to the family
$\{l_p\,|\,p\in\left[1,\infty\right[\}\cup\{c_0\}$. Then
\begin{xxalignat}{2}
&\qquad\qquad\qquad 
K_0\bigg(\allop\bigg(\bigoplus_{j=1}^{n+1} \spx_j\bigg)\bigg) =
\{0\} \qquad \& \qquad K_1\bigg(\allop\bigg(\bigoplus_{j=1}^{n+1}
\spx_j\bigg)\bigg) \cong \inte^n. & \Box
\end{xxalignat}
\end{cor}

As we have seen, every pair of torsion-free, finitely generated, 
commutative groups can arise as
$(K_0(\allop(\spx)),K_1(\allop(\spx)))$. Our last result shows that
$K_0(\allop(\spx))$ may have torsion.

\begin{prop}\label{gm44}
For every $k\in\{2,3, \ldotscor \}$, there is a Banach space $\spx$ for
which $K_0(\allop(\spx))$ contains an element of order $k$.
\end{prop}
\beginpf
In \cite[\S 4.4]{gm2}, Gowers and Maurey describe how to construct a
Banach space $\spx$ satisfying: $\spx^m$ is linearly homeomorphic to
$\spx^n$ if and only if $m \equiv n \pmod{k}$. Hence, it follows from
\eqref{thirdeq} and Proposition~\ref{equi1} that, for every
$\nu\in\nat$, $\nu\,[I]_0 = 0$ if and only if
$k$ divides $\nu$. Consequently, $[I]_0$ has order $k$.
\eopf 
\bigskip 

\noindent 
\begin{center}
\textbf{Acknowledgements}
\end{center}
I want to thank N.\ J.\ Nielsen for sharing with me his insight in the
theory of Banach spaces and, in particular, for introducing me to the
work of Gowers and Maurey. 

I should also like to thank M.\ R{\o}rdam for
teaching me all I know about $K$-theory, for many profitable discussions
and, most of all, for suggesting to me to compute the
$K$-groups of the ideal $\sing$ of strictly singular
operators. This, in fact, started the whole project.

Part of this work was carried out while I visited the University of
Leeds, England, from August to December 1996. I want to thank
everybody in Leeds for an excellent stay. In particular, I should like
to thank H.\ G.\ Dales for his help with the Banach algebraic aspects
of my work, especially for drawing my attention to the question of
when the short exact sequence $\Sigma_{\spx}$ of \S\ref{splitting}
splits.

Finally, I want to thank the three above-mentioned persons for their
careful reading of preliminary versions of this paper.


\begin{thebibliography}{AEO}
\bibitem[AB]{ab} D.\ Alspach \& Y.\ Benyamini, Primariness of Spaces
of Continuous Functions on Ordinals, \emph{Israel J.\ Math.\ }{\bf 27}
(1977), pp.\ 64--92.
\bibitem[AEO]{aeo} D.\ Alspach, P.\ Enflo \& E.\ Odell, On the
Structure of Separable $\mathcal{L}_p$ Spaces $(1 < p < \infty)$,
\emph{Studia Math.\ }{\bf 60} (1977), pp.\ 79--90.
\bibitem[Bl]{bl} B.\ Blackadar, $K$\emph{-Theory for Operator
Algebras,} M.\ S.\ R.\ I.\ Publications, Springer Verlag, 1986. 
\bibitem[Cu]{cuntz} J.\ Cuntz, $K$-Theory for Certain $C^*$-Algebras,
\emph{Ann.\ of Math.\ }{\bf 113} (1981), pp.\ 181--197. 
\bibitem[EW]{ew} I.\ S.\ Edelstein \& P.\ Wojtaszczyk, On Projections
and Unconditional Bases in Direct Sums of Banach Spaces, \emph{Studia
Math.\ }{\bf 56} (1976), pp.\ 263--276.
\bibitem[Fe]{fe} V.\ Ferenczi, Operators on Subspaces of Hereditarily
Indecomposable Banach Spaces, \emph{Bull.\ London Math.\ Soc.,} to appear.
\bibitem[GM1]{gm1} W.\ T.\ Gowers \& B.\ Maurey, The Unconditional
Basic Sequence Problem, \emph{J.\ Amer.\ Math.\ Soc.\ }{\bf 6} (1993),
pp.\ 851--874. 
\bibitem[GM2]{gm2} W.\ T.\ Gowers \& B.\ Maurey, Banach Spaces with
Small Spaces of Operators, \emph{preprint.}
\bibitem[Ka1]{kato1} T.\ Kato, Perturbation Theory for Nullity
Deficiency and Other Quantities of Linear Operators, \emph{J.\ Analyse
Math.\ }{\bf 6} (1958), pp.\ 273--322.
\bibitem[Ka2]{kato2} T.\ Kato, \emph{Perturbation Theory for Linear
Operators,} Die Grundlehren der mathematischen Wissenschaften in
Einzeldarstellungen {\bf 132}, Springer Verlag, 1966.
\bibitem[LP]{lp} J.\ Lindenstrauss \& A.\ Pe{\l}czy{\a'n}ski, Contributions
to the Theory of the Classical Banach Spaces, \emph{J.\ Funct.\ Anal.\
}{\bf 8} (1971), pp.\ 225--249.
\bibitem[LT1]{lt3} J.\ Lindenstrauss \& L.\ Tzafriri, On the
Complemented Subspaces Problem, \emph{Israel J.\ Math.\ }{\bf 9}
(1971), pp. 263--269.
\bibitem[LT2]{lt1} J.\ Lindenstrauss \& L.\ Tzafriri, \emph{Classical
Banach Spaces I,} Ergebnisse der Mathematik und ihrer Grenzgebiete
{\bf 92}, Springer Verlag, 1977. 
\bibitem[Mau]{m} B.\ Maurey, Sous-espaces compl{\'e}mentes de $L_p$,
d'apr{\`e}s P.\ Enflo, \emph{Seminaire Maurey-Schwartz} 1974--1975,
Ecole Polytechnique.
\bibitem[Mit]{mit} B.\ S.\ Mityagin, The Homotopy Structure of the
Linear Group of a Banach Space, \emph{Russian Mathematical
Surveys} {\bf 25} (1970), pp.\ 59--103. 
\bibitem[Pe{\l}]{pel} A.\ Pe{\l}czy{\a'n}ski, On the Isomorphism
Between $m$ and $M$, \emph{Bull.\ Acad.\ Pol.\ Sci.\ }{\bf 6} (1958),
pp.\ 695--696.
\bibitem[Ro]{ro} H.\ P.\ Rosenthal, On Totally Incomparable Banach
Spaces, \emph{J.\ Funct.\ Anal.\ }{\bf 4} (1969), pp.\ 167--175.
\bibitem[Ru]{ru} W.\ Rudin, \emph{Functional Analysis,} 2nd
edn., McGraw-Hill, 1991.
\bibitem[R{\o}r]{ror} M.\ R{\o}rdam, Classification of Inductive Limits
of Cuntz Algebras, \emph{J.\ reine angew.\ Math.\ }{\bf 440} (1993),
pp.\ 175--200.
\bibitem[Tay]{tay} J.\ L.\ Taylor, Banach Algebras and Topology, in
\emph{Algebras in Analysis} (ed.\ J.\ H.\ Williamson), Academic Press,
1975.
\end{thebibliography}
\end{document}